\documentstyle [12pt,a4] {article}
\input epsf
\newcommand{\eps}{\varepsilon}

\newcommand{\dz}{\stackrel{def}{=}}

\newcommand{\minz}{\min\limits}
\newcommand{\sumz}{\sum\limits}
\newcommand{\lz}{\lambda}
\newcommand{\az}{\alpha}

\newcommand{\tr}{\triangle}
\newtheorem{cor}{Corollary}
\newtheorem{prop}{Proposition}
\newtheorem{lem}{Lemma}
\begin{document}
\begin {center}
{\bf On Optimality Conditions for Multi-objective Problems with a Euclidean Cone of Preferences} \end{center}
\begin{center} A. Y. Golubin\\
National Research University Higher School of Economics\\
B. Trechsvjatitelsky per., 3,
 Moscow, 109028, Russia,\\ e-mail: agolubin@hse.ru

\end {center}

{\small The
paper suggests a new  --- to the best of the author's knowledge --- characterization of  decisions which are optimal in the multi-objective optimization problem with respect to a 
definite proper preference cone, a Euclidean cone with a prescribed angular radius.  The main idea is to use the angle distances between the unit vector and points of utility space. 
A  necessary and sufficient condition  for the optimality in the form of an equation is derived. The first-order necessary  optimality conditions 
 are also obtained.}\\

{\it Keywords:} Pareto Optimality;  Cone of Preferences; Scalarization
\\[10pt]
 \paragraph*{1. Introduction }
Optimization problems with several objective functions conflicting with one another are encountered in many situations in practice.
In analyzing such a problem, the concept of Pareto-optimal decisions, that cannot be improved for each criterion without deteriorating the
others, plays an important role.
The Pareto optimality notion is used in solving some engineering and finance problems (Steuer \cite{S}),
insurance theory problems (e.g., Golubin \cite{G2}), ets. There is a
large variety of methods for determination of Pareto-optimal
solutions and their generalizations, in which the optimality is understood with respect to various kinds of preference cones, e.g., Miettinen \cite{1}, Jahn \cite{J}. Many of them
are based on a scalarization approach (Miettinen, Branke et al., Nikulin et al. \cite{1,6,N})  that transforms
the initial problem into a single-objective optimization problem.
Usually it involves some parameters that are changed in order to
detect different Pareto-optimal points: positive weights in the linear
scalarization function or, more generally, a composition of the vector objective function and a linear functional from the dual cone of the preference (ordering) cone; norm parameter 
for $L_p$-scalarization in Nikulin et al.  \cite{N}. Another group of methods for approximating the Pareto
frontier for various decision problems with a small number of
objectives (mainly, two) are provided in Ruzika and  Wiecek \cite{32}. Makela et al. \cite{M} investigate different types of zero-order geometric conditions for characterization of trade-off curves.

The present paper suggests a new "angle distance" scalarization technique for the multi-objective problem with a definite kind of preference cones,  so-called Euclidean cones.
On this basis, a necessary and
sufficient condition for the optimality is derived in the
form of an equation without involving any extra decision maker's parameters. The first-order necessary conditions for weak and locally weak optimal points are also derived.

Formally, the multi-objective optimization problem can be written as:
$$F(x)\equiv (F_1(x),\dots,F_n(x)) \to \max\quad \mbox{s.t. }x\in X,\eqno(1)$$
where $X$ is a decision set or a set of admissible points, $F_i(x)$
are scalar objective functions or utilities defined on $X$. Remark that $F_i(x)$ are not supposed to be concave.

A generalization of the well known Pareto optimality notion is the following (see, e.g., Boyd and Vandenberghe \cite[174]{10}):
Let $K$ be a proper cone: it is convex, closed, the interior $Int\, K\ne \emptyset$, $K$ is pointed, i.e., if $x\in K$ and $-x\in K$ then $x=0$.
 A  point $x^*\in
X$ is called optimal with respect to the cone of preferences $K$ (or $K$-optimal) if there is no other  $x\in X$
such that  
$$F(x)\ne F(x^*) \mbox{ and } F(x)- F(x^*)\in K. \eqno (2)$$ 
A  point $x^*\in X$ is called weak $K$-optimal if  there is no other  $x\in X$  such that 
$$F(x)- F(x^*)\in Int\,K.\eqno(3)$$
Our goal is to find necessary and/or sufficient conditions for optimality in problem (1) for a concrete kind of the cone $K$ which is introduced below.
 \paragraph*{2. Model description }
The concept of Pareto optimality has its root in
economic equilibrium and welfare theory. In the economic terms we try to explain specifics of the suggested modification of the Pareto optimality notion.
Given a set $X$ of alternative allocations of
goods or income for a set of $n$ individuals or members of a community, 
a change from one allocation $x$ to another $y$ is reckoned as an "ideal" improvement if each member increases his/her own  utility by  the same quantity. This means that the increment vector $F(y)-F(x)$ 
lies on the half-line $L^I$ generated by the unit vector $(1,\dots,1)^T\in R^n$ (further we will use its normalized variant $r\dz (1/\sqrt{n},\dots, 1/\sqrt{n})^T$).   A change is considered an 
improvement if a measure of discrepancy between $F(y)-F(x)$  and the "ideal" improvement is not greater than a prescribed constant $a$. An allocation is $K$-optimal
when no further improvement can be made.

 The measure of discrepancy, which defines the very cone $K$ of preferences, is proposed to be the following. Let $p_1$ and $p_2$ be, correspondingly, 
the orthogonal projections of $F(y)-F(x)$ on the "ideal equality" half-line $L^I$ and on the hyper-plane orthogonal to the vector $r$. The measure of discrepancy is 
the norm of $p_2$ per unit of the norm of $p_1$, that is, $\|p_2\|/\|p_1\|$, where $\|z\|$ is the Euclidean norm in $R^n$, $\|z\|=\sqrt{\sum_1^n z_i^2}.$ Passing to the angle distance, 
in our case we have that the above-mentioned discrepancy constraint is $\tan (F(y)-F(x),r)= \sqrt{1-\cos^2(F(y)-F(x),r)}/\cos(F(y)-F(x),r)\le a$. Recall, the cosine of the angle 
between non-zero vectors $x$ and $y$  is defined as $$\cos(x,y)= \frac{<x,y>}{\|x\|\|y\|},$$
where $<x,y>$ denotes the scalar product, $<x,y>=\sum_1^n x_iy_i.$ 
In terms of cosine the latter inequality is expressed as
$\cos(F(y)-F(x),r)\ge s$, where $s=1/\sqrt{a^2+1}$. Now we define the preference cone $K$ as 
$$K(s)\dz \{x\in R^n:\,\cos(x,r)\ge s\}\cup\{0\}\eqno(4)$$
under a given $s\in(0,1)$. Thus, $K(s)$ is a set of 
vectors $x$ such that the angle between the "ideal" direction $r$ and 
each $x$ is not greater than $\arccos s$ (see Fig. 1). Such a cone is called  in Boyd and Vandenberghe \cite[449]{10} a Euclidean cone with the axis $r$ and angular radius $\arccos s$. 
\begin{figure}[h]
\centering
\epsfxsize=0.8\hsize\epsfbox{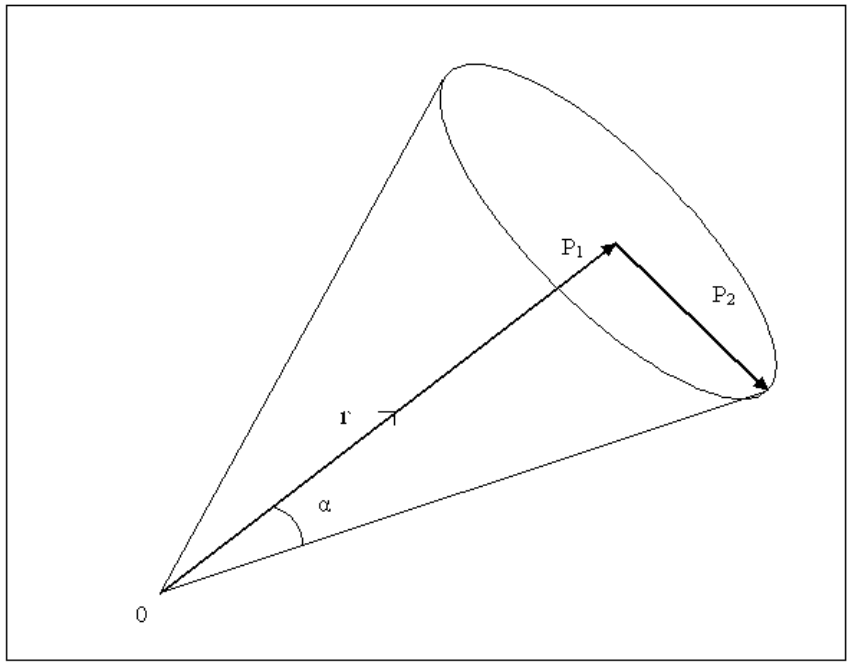} \caption{A preference cone $K(s)$  with angular radius $\az=\arccos s.$}
\end{figure}

 Next we impose reasonable lower and upper boundaries on values of $s$ 
--- it is the same as imposing upper and lower boundaries on the discrepancy limit $a$. The cone of the largest angular radius is supposed  to include the $n$ orts $(0,\dots,0,1,0,\dots,0)^T$ 
of the Pareto preference cone $R^n_+$ as boundary points. Hence, the cosine of the angle between $r$ and each ort is  $1/\sqrt{n}=s$. The cone of the smallest angular radius is supposed  
to include  $n$ orthogonal projections of $r$ on coordinate planes $(x_1,\dots,x_{i-1},0,x_{i+1},\dots,x_n)$ as boundary points. Therefore, the cosine of the angle between $r$ 
and each projection is  $\sum_2^n 1/\sqrt{n}=\sqrt{n-1}/\sqrt{n}=s$. Roughly saying, this cone is "inscribed" of  $R^n_+$, while the previous cone is "circumscribed" of $R^n_+$. 
Then a collection of the preference cones to be considered is
$$\{K(s),\,s\in S\dz [1/\sqrt{n},\sqrt{n-1}/\sqrt{n}]\},$$ 
where $K(s)$ is defined by (4).
It is worth noting that if $n=2$ (bi-objective optimization problem (1)) then the set $\{K(s),\,s\in S\}$ consists of a single cone $R^2_+$ as the interval $S$ converts 
into a point $1/\sqrt{2}$, so the $K(s)$-optimality coincides with the Pareto optimality notion. This particular case was studied in Golubin \cite{G3}.

Denote by $K_U$ and $K_L$ the biggest cone $K(1/\sqrt{n})$ and the smallest cone $K(\sqrt{n-1}/\sqrt{n})$ correspondingly. We call a point $x^*$  {\it upper (lower) optimal} 
if $x^*$ is optimal in (1) with respect to $K_U$ ($K_L$).

From the definition, it is easily seen that any non-zero
$x\in K_L$ necessarily has all non-negative components among which at most one is zero. Indeed, let $x=(0,\dots,0,x_{i+1},\dots,x_n)^T$ and 
$r^1=(0,\dots,0,1/\sqrt{n},\dots, 1/\sqrt{n})^T$ have $i\,(\ge 2)$ zero components. Then $<x,r>=<x,r^1>\le \sqrt{n-i}/(\sqrt{n}\|x\|)<\sqrt{n-1}/(\sqrt{n}\|x\|).$
At the same time, for $n>2$ a vector $x\in K_U$ may have some (not all) negative components. For example,   $x=(-(n-2)/2,1,\dots,1)^T$ is a boundary point of $K_U$ since $\cos(x,r)=1/\sqrt{n}$. 

Returning to the description of $K(s)$-optimality in economic terms, one can say: a $K_L$-improvement of an allocation is that necessarily making at least $n-1$ 
members of a community better off without making the other member worse off, while a $K_U$-improvement may involve  decreases in the utilities of some members, 
however  increases in the utilities of the others make the situation better from the view-point of the community as a whole.


Let $X^*_U,X^*_L,$ and $X^*(s)$ denote, respectively, the sets of all upper optimal, lower optimal, and $K(s)$-optimal points. By construction, 
$K(s_1)\subset K(s_2)$ for $s_1>s_2$, $s_i\in S$. Then, according to the $K(s)$-optimality definition, 
$$X^*_U\subseteq X^*(s_2)\subseteq X^*(s_1)\subseteq X^*_L. \eqno(5)$$
As compared with the Pareto optimality notion, where the preference cone is $R^n_+$, the cone $K_L\subset R^n_+$ and $R^n_+\subset K_U$. 
It leads to that $X^*_U\subseteq X^*_{PO}$ and $X^*_{PO}\subseteq X^*_L $.

Further, when deriving optimality conditions, we will need not only the cone $K(s)$ but also the dual cone of it. Recall that the dual cone of a cone $K$ is 
the set $K^*=\{x\in R^n:\,<x,y> \ge 0 \mbox{ for all } y\in K\}$. Below we give a simple description of the dual of a Euclidean cone with an arbitrary axis.
\begin{lem}
Let $s\in(0,1)$, $q\in R^n$ with $\|q\|=1$, and  $K(q,s)\dz \{x\in R^n:\,\cos(x,q)\ge s\}\cup\{0\}.$ Then the dual cone $$K^*(q,s)=K(q,\sqrt{1-s^2}).\eqno(6)$$ 
\end{lem}
{\bf Proof.} 
By the definition of the cone, it suffices to consider only the vectors of $K(q,s)$ and $K^*(q,s)$ that have a unit norm. Thus, we focus on describing the set 
$E=\{x\in R^n:\,\|x\|=1, <x,y>\ge 0 \mbox{ for all } y\in K(q,s) \mbox{ such that } \|y\|=1\}$. 
First, study the two-dimension case $n=2$. Geometrically, it is clear  that $x\in E$ as long as the angle between $x$ and "the worst" vector $y*\in K(q,s)$ (with a unit length), 
which has the largest angle distance from $x$, is not greater than $\pi/2$. The cosine of an angle $\az$ between $y^*$ and $q$ is $\cos\az=s$, 
therefore\linebreak  $\cos(x,q)=<x,q>=\ge\cos(\pi/2-\az)=\sqrt{1-s^2}$. Thus, (6) is true for $n=2$.

Now proceed with the case $n> 2$. Note, first of all, that $q\in E$ and $-q\notin E$. Fix any $x\in E$, $x\ne q$, and prove that $<x,q>\ge \sqrt{1-s^2}$. 
Let $\Pi$ be a two-dimension plane passing through the vectors $x,q,$ and $0$, i.e. the intersection of all hyper-planes in $R^n$ containing these three points --- note 
that $x$ and $q$ are linearly independent. Denote by $y^*$ any minimum point in the problem
$$\min\, <x,y>\, \mbox{ s.t. } <y,q>\ge s, \|y\|=1 \eqno(7)$$
and show that $y^*\in \Pi$. Consider, at first, an auxiliary problem
with a wider and convex set of admissible points:
$$\min\, <x,y>\, \mbox{ s.t. } <y,q>\ge s, \|y\|\le 1 \eqno(8)$$
Due to convexity of the set $\{y:\, \|y\|\le 1\}$, linearity of both the goal function $<y,q>$ and inequality $<y,q>\ge s $, we have (see, e.g., Bazaraa and Shetty \cite{B}): If $y'$ solves (8), 
then there exists $\lz\in [0,\infty)$ such that $y'$ solves the problem $\minz_{\|y\|\le 1} <x,y>-\lz<q,y>$. Whence
$y'=-(x-\lz q)/\| x-\lz q\|$, therefore $y'$ is also a solution to (7). Thus, $y^*$ is a linear combination of vectors $x$ and $q$, hence, $y^*$ belongs to $\Pi$. 
The latter brings us to the two-dimension case considered above, so $<x,q>\ge \sqrt{1-s^2}$.\\
By analogous reasonings, it can easily be shown that if $x\notin E$, $x\ne -q$, then $<x,q> < \sqrt{1-s^2}$.\\
To sum up, a vector $x$ (of a unit norm) belongs to $E$ if and only if (iff) $<x,q>\ge \sqrt{1-s^2}$, which completes the proof. $\Box$

By definition, $K_U=K(1/\sqrt{n})$ and $K_L=K(\sqrt{n-1}/\sqrt{n})$. Therefore, Lemma 1 gives, in particular, that $K_U$ and $K_L$ are dual cones, 
$$K_U^*=K_L \mbox{ and } K_L^*=K_U.$$
Also, the only self-dual cone of the collection $\{K(s),\,s\in S\}$ is
$K(1/\sqrt{2})$. This is a Lorentz cone (Dattorro \cite[p. 92]{11}) with the axis $r$, its aperture or, in other words, double angular radius is a right angle.
 \paragraph*{3. An "angle distance" scalarization of the problem }
 Let us reformulate the above-given definition (3) of weak $K$-optimality for the case $K=K(s)$: for any $x\in X$ such that $F(x)\ne F(x^*)$ the cosine  
of the angle between $F(x)-F(x^*)$ and $r$ is not greater than $s$ (see Fig. 2), that is,
$$\frac{\sum_1^n( F_i(x) - F_i(x^*))/\sqrt{n}}{\|F(x)-F(x^*)\|}\le s.\eqno(9)$$
\begin{figure}[h]
\centering
\epsfxsize=1.0\hsize\epsfbox{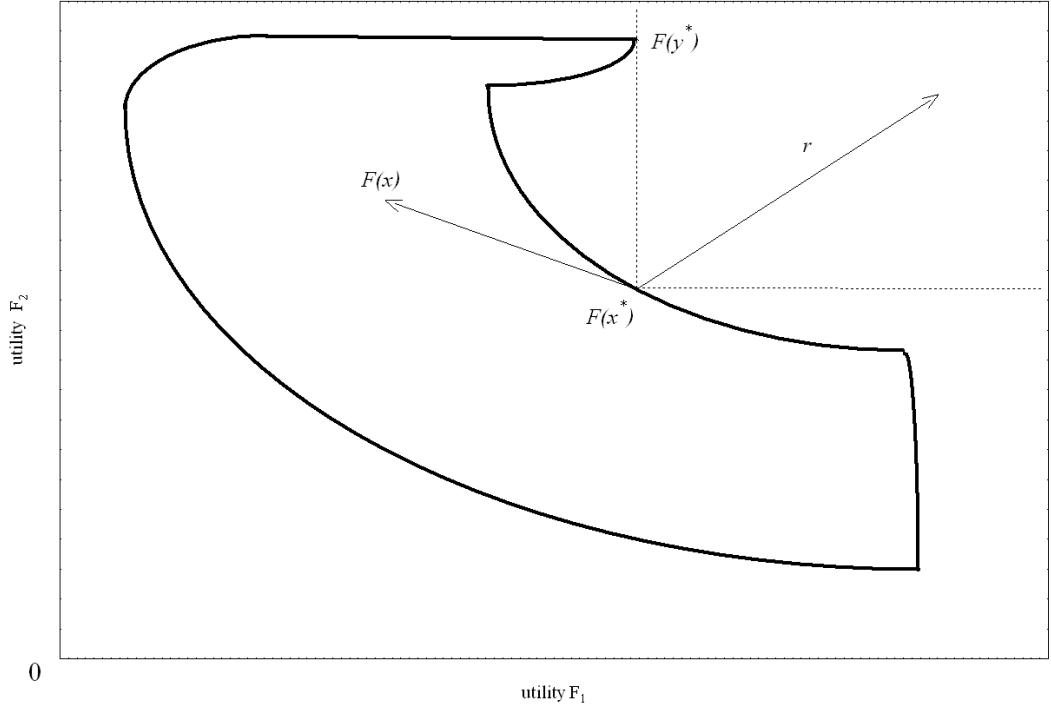} \caption{A non-convex
utility space in $R^2$ with a weak optimal vector $F(x^*).$ }
\end{figure}

Define a scalar function on $X$,
$$G(x)=\sup\limits_{y\in X}\sumz_{i=1}^n(F_i(y) - F_i(x))-s\sqrt{n}\|F(y)-F(x)\|.\eqno(10)$$
By construction, $G(x)\ge 0$ for any $x\in X$ and takes values in the extended real half-line $R_+\cup \{\infty\}$.
Now the necessary and sufficient condition (9) for weak $K(s)$-optimality of $x^*$ can be rewritten as 
$G(x^*)\le 0$.  
Taking into account that supremum in the right-hand side of (10) is attained, in particular, at $y=x^*$,  the latter 
inequality is equivalent to $G(x^*)= 0$.  Thus, we have proved the following proposition

\begin{prop}
A point $x^*$ is weak $K(s)$-optimal iff  $x^*$ is a root of the equation
$$G(x)=0,\eqno(11)$$
where $G(x)$ is defined in \rm{(10)}.
\end{prop}
To find the $K(s)$-optimal (strong) solutions, return to condition (2).  This means that the cone 
$K(s)+F(x^*)$ has no common point with the utility space ${\cal F}\dz \{F(x):\,x\in X\}$, except  $ F(x^*)$.
So, all we need is to find a  weak Pareto-optimal point and to exclude the situation like that depicted on Fig. 2.
\begin{prop}
A point $x^*$ is $K(s)$-optimal iff  $x^*$ is a root of  \rm{(11)}, and
maximum in the problem
$$\max\limits_{y\in X}\sumz_{i=1}^n(F_i(y) - F_i(x))-s\sqrt{n}\|F(y)-F(x)\|\| \eqno (12)$$
is attained at a "unique" point in the sense that if $y^*$ gives
maximum in \rm{(12)} then $F(y^*)=F(x^*)$.
\end{prop}

{\bf Remark 1.}\\
One can easily verify that the role of the "angle distance" scalarization  introduced in (10) and providing
necessary and sufficient conditions  for optimality  in the multi-objective optimization problem (with respect to the cone $K(s)$)  can be played in the case of Pareto optimality 
(with the preference cone $R^n_+$) by a maximin scalarization, 
 $G_1(x)=\sup\limits_{y\in X}\minz_{i=1,\dots,n} F_i(y) -F_i(x)$. The equation corresponding to (11) is then  $G_1(x)=0$.
The functional  $G_1(x)$ is different from (10) even in the case 
$n=2$, where $K(s)=R^2_+$, and generally does not seem convenient for  applications because $G_1(x)$  is not smooth. 
\paragraph*{4. Zero-order optimality conditions}

First, we consider conditions for upper optimality, where, recall, the cone of preferences $K_U=K(1/sqrt{n})$ is the biggest cone of the considered family $\{K(s),\,s\in S\}.$ 
The next statement deals with a zero-order condition for  weak optimality of some point $x^*$, i.e.,  a condition for
solvability of equation (11) with respect to $x^*$, where $s=1/\sqrt{n}$. Denote by
$\tr_i^*(x)= F_i(x) - F_i(x^*),\, i=1,\dots,n,$ where $x^*\in X$ and $x\in X$.
\begin{prop}
A point $x^*$ is weak upper optimal iff  for all $x\in X$ such that  $\sumz_{i=1}^n \tr_i^*(x)>0$, if any, it holds
that $\sumz_{i\ne j}\tr_i^*(x) \tr_j^*(x)\le 0.$
\end{prop}
{\bf  Proof.}
A point $y=x^*$ is a maximizer in problem (12), where now $s\sqrt{n}=1$, iff
$$\sqrt{\sumz_{i=1}^n(\tr_i^*(x))^2}\ge \sumz_{i=1}^n\tr_i^*(x)\eqno(13)$$ for any $x\in X$. If $x$ is such that
$\sumz_{i=1}^n\tr_i^*(x)\le 0$ then (13) holds.  If $\sumz_{i=1}^n\tr_i^*(x)> 0$ then, after  squaring both parts of
(13), we have that (13)  holds iff \linebreak  $2\sumz_{i\ne j}\tr_i^*(x) \tr_j^*(x)\le 0.$ $\Box$

An analogous statement with respect to the (strong) upper optimality follows from the fact that in this case (9) converts into the strict inequality.
\begin{prop}
A point $x^*$ is upper optimal iff  for all $x\in X$ such that $F(x)\ne F(x^*)$ and $\sumz_{i=1}^n \tr_i^*(x)\ge 0$, if any, it holds
that $\sumz_{i\ne j}\tr_i^*(x) \tr_j^*(x)< 0.$
\end{prop}

Return to the general case of the preference cone $K(s),\,s\in S$. The proposition below is, actually, an application of the known results to  our case of the concrete kind of a preference cone.
\begin{prop}
Let $\lz\in Int\,K(\sqrt{1-s^2})$ and $ x^*$ be a maximizer in the problem
$$\max \,\sumz_{i=1}^n \lz_i F_i(x)\quad \mbox{s.t. } x\in X.\eqno(14)$$
Then $ x^*$ is $K(s)$-optimal.
\end{prop}
The  proof consists in a repetition of reasonings in Boyd and Vandenberghe \cite[p. 178]{10} and observing that $K^*(s)= K(\sqrt{1-s^2})$ by Lemma 1.

Since $K_U^*=K_L$ and $K_L^*=K_U$, the next statement directly follows from Proposition 5.

\begin{cor}
Let $\lz\in Int\,K_L$ $(\lz\in Int\,K_U)$   and $ x^*$ be a maximizer in  problem {\rm (14)}.
Then $ x^*$ is $K_U$-optimal ($K_L$-optimal).
\end{cor}

{\bf Remark 2.}\\
 According to the definitions of the cones $ K_L$ and $ K_U$, any weight vector $\lz\in Int\,K_L$ is necessarily positive (component-wise), while
$\lz\in Int\,K_U$ may have some negative components. For instance, in the case $n=4$: a vector $\lz=(1,1,2,2)^T\in Int\,K_L$ as 
$\cos(\lz,r)=6/\sqrt{40}> s=\sqrt{n-1}/\sqrt{n}=\sqrt{3}/\sqrt{4}$; a vector $\lz=(-\eps,-\eps,1,1)^T$, where $\eps\in(0,2-\sqrt{3})$, belongs to 
$Int\,K_U$ as $\cos(\lz,r)> s=1/\sqrt{n}=1\sqrt{4}$. 

The existence of an upper optimal point (and, therefore, any $K(s)$-optimal point for $s\in S$ (see (5)) is guarantied by solvability of problem (14) with positive $\lz\in Int\,K_L$ which, in turn, 
is guarantied by compactness of $X$ and upper semi-continuity of all $F_i(x)$.
 
\paragraph*{5. First-order necessary conditions for optimality }
Let $x^*$ be a weak  $K(s)$-optimal point, i.e., a root of (11). Denote by $y^*$ a maximum point in (12). It is 
easily seen that $y^*$ is also weak  $K(s)$-optimal. We will call such a pair $(x^*, y^*)$  a 
{\it weak  $K(s)$-optimal pair}. Of course, $y^*$ can always be taken equal to $x^*$, so  $( x^*, x^*)$  
is always a weak  $K(s)$-optimal pair. A more interesting situation is that where $x^*$ is not a unique solution 
to maximization problem (12).
In the sequel of this section we suppose that the decision set $X\subseteq R^k$ and utility functions
$F_i(x),\,i=1,\dots,n$ are differentiable  on $R^k$.
\begin{prop}
Let $( x^*, y^*)$ be a weak  $K(s)$-optimal pair and $y^*$ be an internal point of $X$.
Then
\setcounter{equation}{14}
\begin{eqnarray}
&&\sumz_{i=1}^n F_i'(y^*)[\sumz_{j=1}^n \tr_j^*(y^*)-s^2n \tr_i^*(y^*)]=0,\\
&& \sumz_{i=1}^n \tr_i^*(y^*) -s\sqrt{n}\|F(y^*)-F(x^*)\|=0.
\end{eqnarray}
\end{prop}
{\bf  Proof.} Suppose, at first,  that $F(y^*)\ne F(x^*)$. Since
$y^*$ solves problem (12), the first-order optimality condition is
$$ \sumz_{i=1}^n F_i'(y^*)- s\sqrt{n}\frac{\sumz_{j=1}^n F_j'(y^*)\tr_j^*(y^*)}{\|F(y^*)-F(x^*)\|} =0,\eqno(17)$$
where, recall, $\tr_i^*(y^*)= F_i(y^*) - F_i(x^*).$ From Proposition
1 it follows that
$$ \sumz_{j=1}^n F_j'(y^*)\tr_j^*(y^*)= s\sqrt{n}\|F(y^*)-F(x^*)\|.\eqno(18)$$ After substituting the expression for $\|F(y^*)-F(x^*)\|$ into (17), we obtain
\begin{eqnarray*}
&&\sumz_{i=1}^n F_i'(y^*)[\sumz_{j=1}^n \tr_j^*(y^*)-s^2n \tr_i^*(y^*)]=0.
\end{eqnarray*}
The latter relation admits the degenerated case $F(y^*)=F(x^*)$ also. Taking (18) into account, we complete the proof. $\Box$

{\bf Remark 3.}\\
The statement of Proposition 6 becomes trivial  if  a maximum point
$y^*$  in (12) corresponds to the same point in utility space as
$x^*$, $F(y^*)=F(x^*)$. Nevertheless, Proposition 6 provides an
informative necessary condition in the case where $x^*$ is just weak
$K(s)$-optimal, but not strong $K(s)$-optimal, as shown on Fig. 2. In the case of upper optimality, where $s=1/\sqrt{n}$, equation (15) becomes simpler
$$\sumz_{i=1}^n F_i'(y^*)\sumz_{j\ne i} \tr_j^*(y^*)=0. \eqno(19)$$
\quad Below we will need  a notion of local $K(s)$-optimum.  A point $x^*$ is called   local weak $K(s)$-optimal if there exists an
$\eps$-neighborhood $O_{\eps}(x^*)$ of this point such that $x^*$ is
weak $K(s)$-optimal with respect to a smaller decision set
$O_{\eps}(x^*)\cap X$. 

Like Proposition 5, the next proposition applies a known
result on optimality conditions for the multi-objective optimization problem to the considered Euclidean cones $K(s)$ of preferences. 
\begin{prop}
Let $ x^*$ be a local weak  $K(s)$-optimal point and  an internal point of $X$.
Then there exists a vector $\lz\in K(\sqrt{1-s^2})\backslash \{0\}$  such that
$$\sumz_{i=1}^n \lz_i F_i'(x^*)=0.\eqno(20)$$
\end{prop}
{\bf  Proof.} Denote by $F'(x^*)$ a matrix $n\times k$, its rows are  gradients $F_i'(x^*)$,  $i=1,\dots,n.$ Firstly, prove that there is no $h\in R^k$ such that $F'(x^*)h\in Int\,K(s)$.
 Suppose the contrary, then $F(x^*+th)-F(x^*)\in Int\,K(s)$ for all sufficiently small $t>0$, which
contradicts to the above-supposed local weak $K(s)$-optimality of
$x^*$. Applying the so-called alternative theorem (see, e.g., Boyd and Vandenberghe \cite [p.54]{10}), we have that there exists $\lz\ne 0$ such that  $\lz\in K^*(s)$ and 
$\lz^T F'(x^*)=0$ (i.e.  $\sumz_{i=1}^n \lz_i F_i'(x^*)=0$).  By Lemma 1, $K^*(s)= K(\sqrt{1-s^2})$.
$\Box$

\begin{cor}
Let $ x^*$ be a local weak  upper optimal (lower optimal) point and  an internal point of $X$.
Then there exists a vector $\lz\in K_L\backslash \{0\}$ $(\lz\in K_U\backslash \{0\})$   such that \rm{(20)} holds.
\end{cor}


 \end{document}